\documentclass[a4j,10pt,onecolumn,oneside,notitlepage,openany,final]{report}
\usepackage{amsmath}
\usepackage{amsthm}
\usepackage{amsfonts}
\usepackage{amssymb}
\usepackage{bm}
\usepackage{graphicx}
\usepackage{longtable}

\newtheorem{thm}[]{Theorem} 
\newtheorem{Def}[thm]{Definition}
\newtheorem{lem}[thm]{Lemma}
\newtheorem{cor}[thm]{Corollary}
\newtheorem{prp}[thm]{Proposition}
\newtheorem{rem}[thm]{Remark}

\begin{document}
\title{\textbf{\mathversion{bold} The punctual Hilbert schemes for the curve singularities of types $E_6$ and $E_8$}}
\author{Yoshiki S\={O}MA,\ Masahiro WATARI}
\date{}
\maketitle
\begin{abstract}
Pfister and Steenbrink studied punctual Hilbert schemes for irreducible  curve singularities. 
In particular, they  investigated the structure of  special punctual Hilbert schemes for certain monomial curve singularities. 
In this paper, we study the punctual Hilbert schemes of all degrees for the curve singularities of types $E_6$ and $E_8$. 
For our analysis, we introduce computational algorithms to decompose a punctual Hilbert schemes into affine cells. 
We also use the theory of Gr\"obner basis and known results about the compactified Jacobian of singular curves to prove our main theorems.
\end{abstract}
\noindent
\textbf{Keywords} punctual Hilbert schemes, curve singularities of types $E_6$ and $E_8$, compactified Jacobians of singular curves\\
\textbf{Mathematics Subject Classification (2000)} 14C05, 14H20

\section{Introduction} 
Let $\mathcal{O}$ be the complete local ring of a irreducible curve singularity over an algebraically closed field $k$ of characteristic 0. 
We denote by $\overline{\mathcal{O}}$ and $\delta$ the normalization of $\mathcal{O}$ and the $\delta$-invariant of $\mathcal{O}$ respectively. 
Pfister and Steenbrink  \cite{PS} defined a special subset $\mathcal{M}$ of the Grassmannian Gr$\left(\delta,\overline{\mathcal{O}}/I(2\delta)\right)$ where $I(2\delta)$ is the set of all elements in $\mathcal{O}$ 
whose orders are greater than or equal to $2\delta$. 
It is known that $\mathcal{M}$ is  a projective variety defined by Pl\"{u}cker relations and additional linear equations (see \cite{PS}). 
We call it the Pfister-Steenbrink variety (PS variety) for a given singularity. 
By using its intersection with Schubert cells, they investigated the structure of $\mathcal{M}$ for the curve singularities with monomial semigroups. 
The punctual Hilbert scheme $\mathcal{M}_r$ of degree $r$ was also constructed  as a connected component of $\mathcal{M}$. 
It is a projective variety which parametrizes the ideals of codimension $r$ in $\mathcal{O}$.  
The PS variety  $\mathcal{M}$ coincides with punctual Hilbert schemes whose degrees are greater than or equal to 2$\delta$ (Corollary\,\ref{range}).

In the present paper, we study the structure of all punctual Hilbert schemes for  the curve singularities of types $E_6$ and $E_8$
(i.e. the  curve singularities whose local rings are $k[[t^3,t^4]]$ and $k[[t^3,t^5]]$ respectively).  
The PS varieties for these singularities were originally studied in \cite{PS}. 
See also \cite{S} which was  the preliminary version of this paper.
Our main theorems are stated as follows:
\begin{thm}\label{main 1}
For the curve singularity of type $E_6$, 
its punctual Hilbert schemes are given by the following table$:$
\begin{center}
\begin{tabular}{c||c|c|c|c|c|c|}
$r$&$1$&$2$&$3$&$4$&$5$&$\ge 6$\\
\hline
$\mathcal{M}_r$&$\mathbb{P}^0$&$\mathbb{P}^1$&$\mathbb{P}^2$&$\mathbb{P}^2\cup X_1$&$\mathbb{P}^2\cup \mathbb{P}^2$&$X_2$\\
\hline
$\mathrm{Sing}(\mathcal{M}_r)$&\multicolumn{3}{|c|}{$\emptyset$}&\multicolumn{3}{|c|}{$\mathbb{P}^1$}\\
\multicolumn{7}{c}{\textbf {Table\,}{\mathversion{bold}$1$.}}
\end{tabular}
\end{center}
The variety $X_1$ $($resp. $X_2)$ in Table\,$1$ is a rational projective  surface $($resp. a rational projective  threefold$)$. 
Their defining equations are listed in Section\,$\ref{E6}$. 
\end{thm}

\begin{thm}\label{main 2}
For the curve singularity of type $E_8$, 
its punctual Hilbert schemes are given by the following table$:$

\begin{center}
\begin{tabular}{c||c|c|c|c|c|c|c|c|}
$r$&$1$&$2$&$3$&$4$&$5$&$6$&$7$&$\ge 8$\\
\hline
$\mathcal{M}_r$&$\mathbb{P}^0$&$\mathbb{P}^1$&$\mathbb{P}^2$&$\mathbb{P}^2\cup X_3$&$\mathbb{P}^2\cup \mathbb{P}^2\cup X_4$&$X_5\cup X_6$&$X_7\cup X_8$&$X_9$\\
\hline
$\mathrm{Sing}(\mathcal{M}_r)$&\multicolumn{3}{|c|}{$\emptyset$}&$\mathbb{P}^1$&$\mathbb{P}^1\cup \mathbb{P}^1$&\multicolumn{3}{|c|}{$\mathbb{P}^2\cup \mathbb{P}^2$}\\
\multicolumn{9}{c}{\textbf {Table\,}{\mathversion{bold}$2$.}}
\end{tabular}
\end{center}

The all  varieties $X_i$ $(i=3,\cdots,8)$ in Table\,$2$ are rational and projective.
Their dimensions are given by 
\begin{equation*}
\dim X_i=
\begin{cases}
2\quad\text{for $i=3$,}\\
3\quad\text{for $i=4,\ldots,8$,}\\
4\quad\text{for $i=9$.}\\
\end{cases}
\end{equation*}

\end{thm} 

The PS variety for a curve singularity was studied in another point of view. 
Rego \cite{R} introduced the compactified Jacobian of singular curves. 
He also constructed the Jacobi factor for a curve singularity. 
For a given singularity, the Jacobi factor coincides with the PS variety in his construction. 
Beauville \cite{B} proved that the Euler numbers of the  Jacobi factors for the irreducible curve singularities with Puiseux characteristic $(p,q)$ are given by $\frac{1}{p+q}\binom{p+q}{p}$. 
Piontkowski  \cite{P}  computed the Euler numbers of the Jacobi factors for  irreducible curve singularities with  Puiseux characteristics $(4,2q,s)$ and $(6,8,s)$. 
We will show that the Beauville's theorem for the cases $(p,q)=(3,4)$ and $(3,5)$ (i.e. the case for curve singularities of types $E_6$ and $E_8$) also follows from our arguments. 

The paper is laid out as follows: 
In Section\,\ref{PS-theory}, we briefly recall the Pfister-Steenbrink theory for punctual Hilbert schemes of curve singularities. 
We also fix notations and prove some lemmas needed later. 
In Section\,3, we introduce computational algorithms to determine the affine cells of $\mathcal{M}_r$. 
Some useful results of Gr\"obner basis are also listed. 
Finally, by using them, we prove Theorem\,\ref{main 1} and \ref{main 2}  in Section\ \ref{E6} and \ref{E8} respectively. 
Some known results of  compactified Jacobians of singular curves are also used to prove Theorem\,\ref{main 2}. 
\vspace{3mm}

\noindent
\textbf{Acknowledgement} 
The authors would like to express his sincere gratitude to Professor Fumio Sakai for his
valuable advices and warm encouragement during the preparation of the present article.
 They would also like to thank Professor Matsuda for his useful advices.

\section{Preliminaries}\label{PS-theory}
In this paper, we only consider monomial curve singularities defined below. 
\begin{Def}\label{Def of monomial}
A monomial curve singularity is an irreducible curve singularity whose local ring is isomorphic to $k[[t^{a_1},\ldots,t^{a_m}]]$ for some $a_1,\ldots,a_m\in \mathbb{N}$. 
\end{Def}
\begin{rem}
Without loss of generality, we may assume that {\rm gcd}$(a_1,\ldots,a_m)=1$ in Definition\,$\ref{Def of monomial}$. 
\end{rem}

Let $\mathcal{O}=k[[t^{a_1},\ldots,t^{a_m}]]$ be the local ring of  a monomial curve singularity. 
Note that its normalization $\overline{\mathcal{O}}$ is $k[[t]]$. 
We call $\Gamma:=\{\text{ord}(f)\,|\,f\in \mathcal{O}\}$ the \emph{semigroup} of $\mathcal{O}$.  
A positive integer $\delta:=\dim _{k}(\overline{\mathcal{O}}/\mathcal{O})$ is called the \emph{$\delta$-invariant} of $\mathcal{O}$.  
For a natural number $n$, set $\overline{I}(n):=\{\,f\in \overline{\mathcal{O}}|\, \text{ord}(f)\ge n\}$ and $I(n):=\overline{I}(n)\cap \mathcal{O}$. 
Setting ord$(0)=\infty$, we regard $\overline{I}(n)$ (resp. $I(n)$) as an ideal of $\overline{\mathcal{O}}$ (resp. $\mathcal{O}$).    
Define the \emph{conductor} $c$ of $\mathcal{O}$ to be $\min\{n|\,\overline{I}(n)\subset \mathcal{O}\}$. 
It is known that  $\delta+1\le c\le 2\delta$ and $c=2\delta$ if and only if $\mathcal{O}$ is Gorenstein (cf. \cite{Serre}).  
For an ideal $I$ of $\mathcal{O}$, we call  $\Gamma(I):=\{\text{ord}(f)|\,f\in \mathcal{O}\}$ the \emph{order set} of $I$.

A subset $\Delta\subset \mathbb{Z}_{\geq 0}$ is called a $\Gamma$-semi-module, if $\Delta+\Gamma\subset \Delta$.
If a $\Gamma$-semi-module $\Delta$ is minimally generated by $\alpha_1,\cdots,\alpha_m$ (i.e $\Delta=\sum_{i=1}^m(\alpha_i+\Gamma)$ 
and $\Delta\supsetneq \sum_{i=1,i\neq j}^m(\alpha_i+\Gamma)$ for $\forall j\in\{1,\ldots,m\}$), then we write $\Delta=\langle \alpha_1,\cdots,\alpha_m\rangle_{\Gamma}$. 
We denote by $\mathcal{I}(\Delta)$ the set of all ideals of $\mathcal{O}$ whose order sets are  $\Delta$. 
Note that $\mathcal{I}(\Delta)\neq \emptyset$ if and only if $\Delta\subset \Gamma$. 
The following lemma is obvious:

\begin{lem}\label{normal form}
Let $\Delta=\langle \alpha_1,\cdots,\alpha_m\rangle_{\Gamma}$ be a $\Gamma$-semi-module. 
If $\mathcal{I}(\Delta)\neq \emptyset$, then  any ideal  in $\mathcal{I}(\Delta)$ is minimally generated by the  elements of the form
\begin{equation*}
f_i:=t^{\alpha_i}+\sum_{j\in  \Gamma\setminus \Delta,j>\alpha_i}a_{i,j}t^j\quad (a_{i,j}\in k,\ i=1,\ldots,m). 
\end{equation*}
\end{lem}

\noindent
For a positive integer $r$, set
$\mathcal{I}_r:=\{I|\,\text{$I$ is an ideal of $\mathcal{O}$ with $\dim_k\mathcal{O}/I=r$}\}$. 
\begin{lem}\label{colength}
An ideal in $\mathcal{O}$ belongs to  $\mathcal{I}_r$ if and only if we have $\sharp \{\Gamma\setminus \Gamma(I)\}=r$. 
\end{lem}
\begin{proof}
It is clear that $I$ belongs to $\mathcal{I}_{r}$ if and only if 
$$\mathcal{O}/I=\left\{\overline{a_0+a_1t^{d_1}+\cdots+a_{r-1}t^{d_{r-1}}}\big |\,a_i\in k, d_i\in \Gamma\setminus \Delta, d_1<\cdots<d_{r-1}\right\}$$ holds. 
Thus, we have  $\sharp \{\Gamma\setminus \Gamma(I) \}=r$. 
\end{proof}

\noindent
If $\Delta$ is a $\Gamma$-semi-module such that $\mathcal{I}(\Delta)\neq \emptyset$, then 
all ideals in $\mathcal{I}(\Delta)$ have same codimension by  Lemma\,\ref{colength}. 
So the following fact folds:
\begin{lem}\label{colength2}
Let $\Delta$ be a $\Gamma$-semi-module such that $\mathcal{I}(\Delta)\neq \emptyset$. 
A set $\mathcal{I}(\Delta)$  is contained in  $\mathcal{I}_r$ if and only if we have $\sharp \{\Gamma\setminus \Delta\}=r$. 
\end{lem}

\begin{prp}\label{decomposition lem}
There exists a finite number of distinct $\Gamma$-semi-modules $\Delta_{r,1},\cdots,\Delta_{r,l_r}$ such that
\begin{equation}\label{decomposition}
\mathcal{I}_r=\bigcup_{i=1}^{l_r} \mathcal{I}(\Delta_{r,i}).
\end{equation}
\end{prp}
\begin{proof}
The finiteness of the number of  $\Gamma$-semi-modules holds trivially, as there exists only a finite number of semigroups in $\mathbb{N}$ of fixed colength. 
It is clear that (\ref{decomposition}) is a disjoint union.
\end{proof}
\begin{rem}
By Lemma\,$\ref{colength2}$, the $\Gamma$-semi-modules $\Delta_{r,i}'s$ in  $(\ref{decomposition})$ are an invariant for  codimension.
\end{rem}

Let $\rm{Gr}\left(\delta,\overline{\mathcal{O}}/I(2\delta)\right)$ be the Grassmannian which consists of $\delta$-dimensional linear subspaces of $\overline{\mathcal{O}}/I(2\delta)$.
For $V\in\rm{Gr}\left(\delta,\overline{\mathcal{O}}/I(2\delta)\right)$, define a multiplication by $\mathcal{O}\times V\ni (f,\overline{v})\mapsto \overline{fv}\in V$. 
Set
$$
\mathcal{M}:=\left\{V\in \rm{Gr}\left(\delta,\overline{\mathcal{O}}/I(2\delta)\right)\big|\,\text{$V$ is an $\mathcal{O}$-submodule w.r.t the multiplication}\right\}.
$$
Consider a composition map
$$
\psi:\,\mathcal{M}\ \rightarrow \text{Gr}(\delta,2\delta)\rightarrow \text{M}_{\delta,2\delta}(k)/\sim \ \rightarrow \mathbb{P}^{N
}$$ 
where $\mathrm{Gr}(\delta, 2\delta)$ is the Grassmannian which consists of $\delta$-dimensional linear subspaces of $k^{2\delta}$, $M_{\delta,2\delta}(k)$ is the set of all $\delta\times 2\delta$ matrices over $k$ and the equivalence relation $\sim$ is the similarity of matrices. 
For a formal power siries $f=\sum_{j=0}^\infty a_{j}t^j$ in $\overline{\mathcal{O}}$, we denote its coset in $\overline{\mathcal{O}}/I(2\delta)$ by $\overline{f}=\sum_{j=0}^{2\delta-1}a_{j}\tau^j$ where $\tau\equiv t$ mod $I(2\delta)$. 
For a coset $\overline{f}$, define its order $\mathrm{ord} (\overline{f})$ by 
$\mathrm{ord} (f)$ (resp. $\infty$), if $\mathrm{ord}(f)\le 2\delta-1$ (resp. $\overline{f}=\overline{0}$). 
We use the notation  $[ a_1,\cdots, a_n ] _k$ for a $k$-vector space generated by $a_1,\ldots, a_n$. 
Let $V=[ \overline{f}_1,\cdots, \overline{f}_\delta ] _k$ be an element of $\mathcal{M}$ where $\overline{f}_i=\sum_{j=0}^{2\delta-1}a_{i,j}\tau^j$. 
We identify  $\overline{f}_i$ with the point $\bm{a}_i=(a_{i,0},\cdots,a_{i,2\delta-1})$ in $k^{2\delta}$.  
The first map in $\psi$ is defined by this identification. 
Let $A_V$ be the $\delta \times 2\delta$ matrix whose $i$th row is $\bm{a}_i$. 
We  call   it the \emph{representation matrix} of $V$. 
The second map in $\psi$ sends a $k$-vector space $[\bm{a}_1,\cdots,\bm{a}_{\delta}]_{k}$ to   the coset of  $A_{V}$ . 
We may assume that the coset of $A_{V}$ is represented by the reduced row echelon form. 
The third map in $\psi$ is Pl\"{u}cker embedding with $N=\binom{2\delta}{\delta}-1$.

For $r>0$, Pfister and Steenbrink defined a map $\varphi_r:\mathcal{I}_r\rightarrow \mathcal{M}$ by $\varphi _r (I)=t^{-r}I/I(2\delta)$. 
\begin{prp}[\cite{PS}, Theorem\,3]\label{thm3}
The map $\varphi_r$  is injective for any $r$. 
Furthermore, it is bijective for $r\ge 2\delta$. 
The image $(\psi\circ\varphi_r)(\mathcal{I}_r)$  is a Zariski closed set in $\psi(\mathcal{M})$. 
\end{prp}

\noindent
Put $\mathcal{M}_r:=\varphi_r(\mathcal{I}_r)$. 
Since $\psi$ is injective, we identify $\psi(\mathcal{M})$ and $\psi(\mathcal{M}_r)$ with $\mathcal{M}$ and  $\mathcal{M}_r$ respectively. 
\begin{Def}
We call $\mathcal{M}$ and  $\mathcal{M}_r$ the Pfister-Steenbrink variety $($PS variety$)$ and the punctual Hilbert scheme of degree $r$ for a given curve singularity respectively.
\end{Def} 

\noindent
The following fact follows from Proposition\,\ref{thm3}:
\begin{cor}\label{range}
Any punctual Hilbert scheme $\mathcal{M}_r$ with $r\ge 2\delta$ coincides with the PS variety $\mathcal{M}$.
\end{cor}

\noindent
By Corollary\,$\ref{range}$, it is enough to consider $r$ within  $1\le r\le 2\delta$ for the analysis of $\mathcal{M}_r$. 
Setting $\mathcal{M}_{r,i}:=\varphi_r(\mathcal{I}(\Delta_{r,i}))$ for a component $\mathcal{I}(\Delta_{r,i})$ in (\ref{decomposition}), 
we have the following:

\begin{lem}\label{generators for k-vec.}
Let  $\beta_{1},\ldots,\beta_{\delta}$  be the  first $\delta$ elements in $\Delta$. 
As a $k$-vector space, any element in $\mathcal{M}_{r,i}$ is generated by the elements of the form 
\begin{equation}\label{basis}
\overline{g}_i:=\tau^{\beta_i-r}+\sum_{j\in \Gamma\setminus \Delta,\, j>\beta_i}b_{i,j}\tau^{j-r} \quad (i=1,\ldots,\delta).
\end{equation} 
\end{lem}

\begin{proof}
Our assertion follows from the definition of $\varphi_r$ and Lemma\,\ref{normal form}. 
\end{proof}

\noindent
Since $\psi$ is injective, we also identify $\psi (\mathcal{M}_{r,i})$ with $\mathcal{M}_{r,i}$.
Namely, $\mathcal{M}_{r,i}$ is regarded as the subset of  the punctual Hilbert scheme $\mathcal{M}_r$ parametrizing ideals in $\mathcal{I}(\Delta_{r,i})$.
We set $[a,b]:=\{x\in\mathbb{Z}_{\ge 0}|\,a\le x\le b\}$. 
The following fact is known:
\begin{prp}[\cite{PS}, Corollary of Theorem\,11]\label{affine space}
The set $\mathcal{M}_{r,i}$ is isomorphic to 
the affine space $\mathbb{A}^N$ where $N=\sum_{\gamma\in (\Delta_{r,i}-r)\cap[0,2\delta]}\sharp J_{\gamma}$ where $J_{\gamma}:=[\gamma+1,2\delta-1]\setminus \Delta_{r,i}$.
\end{prp}

The affine cell decomposition of $\mathcal{M}_r$  follows from Proposition\,\ref{decomposition lem} and \ref{affine space}.
\begin{cor}\label{decomposition cor}
The punctual Hilbert scheme  $\mathcal{M}_r$ of degree $r$ has the following affine cell decomposition$:$ 
\begin{equation}\label{decomposition 2}
\mathcal{M}_r=\bigcup_{i=1}^{l_r}\mathcal{M}_{r,i}
\end{equation}
\end{cor}

\begin{rem}
The affine cells $\mathcal{M}_{r,i}$ $(i=1,\ldots, l_r)$ in $(\ref{decomposition 2})$ are induced by the canonical flag 
$0\subset V_1\subset \cdots \subset V_{2\delta}=\overline{\mathcal{O}}/(t^{2\delta})$
where $V_i=\overline{I}(2\delta-i)/I(2\delta)$. For details, see $\cite{PS}$. 
\end{rem}

\noindent
Proposition\,\ref{affine space} also  yields the following: 
\begin{cor}\label{rationality}
If $\mathcal{M}_r$ is irreducible, then it is a rational projective variety.
\end{cor}

\section{Computational algorithms}\label{Algorithms}
The aim of this section is to prove Theorem\,\ref{Computational steps}
which determines the affine cell decomposition (\ref{decomposition 2})  of $\mathcal{M}_r$.  
We freely use the notations introduced in the previous section. 
\begin{lem}\label{alg1}
For a  $\Gamma$-semi-module $\Delta$, we obtain the set $A$ of minimal generators of $\Delta$ by the following algorithm$:$

\noindent
\rm{\texttt{INPUT:}}$\Delta$ \\
\rm{\texttt{OUTPUT:}}$A$\\
\rm{\texttt{DEFINE:}}$A:=\emptyset$, $\Delta:=\Delta$\\
\rm{\texttt{WHILE}} $\Delta\neq \emptyset$ \rm{\texttt{DO}}\\
\quad $A:=A\cup\{\min\{\Delta\}\}$\\
\quad $\Delta:=\Delta\setminus \{\min\{\Delta\}+\gamma\,|\,\gamma\in \Gamma\}$
\end{lem}
\begin{proof}
It is trivial. 
So we omit the proof. 
\end{proof}

\begin{lem}\label{prp r r+1}
Let $\Delta=\langle \alpha_1,\cdots,\alpha_m\rangle_{\Gamma}$ be a $\Gamma$-semi-module.  
If  $\mathcal{I}(\Delta)$ is a component of $\mathcal{I}_r$, then 
$\mathcal{I}(\Delta\setminus \{\alpha_i\})$ is  component of $\mathcal{I}_{r+1}$  for each $i\in\{1,\ldots,m\}$. 
Conversely, if  $\mathcal{I}(\Delta)$ is a component of $\mathcal{I}_{r+1}$, then, for each $\alpha_i$ and $\gamma_1:=\min\{\Gamma\setminus \{0\}\}$,  
$\mathcal{I}(\Delta\cup \{\alpha_i-\gamma_1 \})$ is a component of $\mathcal{I}_r$.
\end{lem}

\begin{proof}
Assume that $\mathcal{I}(\Delta)$ is a component of $\mathcal{I}_r$. 
For any $\alpha_i$, it is clear that $\Delta\setminus \{\alpha_i\}$ is also a $\Gamma$-semi-module. 
Since $\sharp (\Gamma\setminus \Delta)=r$ by  Lemma\,\ref{colength}, we have $\sharp (\Gamma\setminus (\Delta\setminus \{\alpha_i\}))=r+1$. 
Hence the set $\mathcal{I}(\Delta\setminus \{\alpha_i\})$ is  component of $\mathcal{I}_{r+1}$.  
Next assume that the set $\mathcal{I}(\Delta)$ is a component of $\mathcal{I}_{r+1}$. 
Now we have $\alpha _i-\gamma_1\notin \Delta$ for any $i$. 
Indeed, if $\alpha_i-\gamma_1\in \Delta$, then there exist $\alpha_j$ and $\gamma$ in $\Gamma$ such that $\alpha _i-\gamma_1=\alpha_j+\gamma$. 
This fact implies that $\alpha_1,\cdots,\alpha_m$ are not  minimal generators of $\Delta$. 
It contradicts the assumption. 
It is clear that $\Delta\cup \{\alpha_i-\gamma_1\}$ is a $\Gamma$-semi-module and 
$\sharp (\Gamma\setminus (\Delta\cup \{\alpha_i-\gamma_1\}))=r$. 
Hence, $\mathcal{I}(\Delta\cup \{\alpha_i-\gamma_1\})$ is a component of $\mathcal{I}_r$ by Lemma\,\ref{colength}.
\end{proof}
For the decomposition\,(\ref{decomposition}) of $\mathcal{I}_r$, set $\mathfrak{D}_r:=\{ \Delta_{r,1},\cdots,\Delta_{r,l_r}\}$. 
We have the following proposition which determines $\mathfrak{D}_r$ from $\mathfrak{D}_{r-1}$: 
\begin{prp}\label{alg2}
We construct $\mathfrak{D}_r$ from $\mathfrak{D}_{r-1}$ in a finite number of steps given by the following algorithm$:$

\noindent
\rm{\texttt{INPUT:}}
$\mathfrak{D}_{r-1}=\{\Delta_{r-1,1}, \cdots,\Delta_{r-1,l_{r-1}}\}$ where $\Delta_{r-1,i}=\langle\alpha_{i,1},\cdots,\alpha_{i,m(i)}\rangle_{\Gamma}$ $(i=1,\ldots,l_{r-1})$ \\
\rm{\texttt{OUTPUT:}} $\mathfrak{D}_r$\\
\rm{\texttt{DEFINE:}} $\mathfrak{D}_{r}:=\emptyset$\\
\rm{\texttt{FOR}} each $i\in\{1\,\ldots,l_{r-1}\}$ and each $j\in\{1\,\ldots,m(i)\}$ \rm{\texttt{DO}}\\
\quad $\Delta:= \Delta_{r-1,i}\setminus \{\alpha_{i,j}\}$\\
\qquad \rm{\texttt{IF}} $\Delta\notin \mathfrak{D}_{r}$ \rm{\texttt{THEN}} $\mathfrak{D}_{r}:=\mathfrak{D}_{r}\cup \{\Delta\}$
\rm{\texttt{ELSE}} do nothing
\end{prp}

\begin{proof}
Our assertion follows from Lemma\,\ref{prp r r+1}. 
\end{proof}
Let $\mathcal{I}(\Delta_{r,i})$ be a component of the decomposition (\ref{decomposition}) of $\mathcal{I}_r$. 
Recall that $\mathcal{M}_{r,i}=\varphi_r (\mathcal{I}(\Delta_{r,i}))$ is a family of $k$-vector spaces of dimension $\delta$. 
Since  all elements in $\mathcal{M}_{r,i}$ has generators of the form (\ref{basis}) as in Lemma\,\ref{generators for k-vec.}, 
 we regard $\mathcal{M}_{r,i}$ itself as $a$ $k$-vector space generated by $\overline{g}_1,\cdots,\overline{g}_\delta$. 
So we  just write $\mathcal{M}_{r,i}=[ \overline{g}_1,\cdots,\overline{g}_\delta]_k$. 
Here the  coefficients $b_{i,j}$'s in ($\ref{basis}$) are treated as variables. 
Since the set $\{(\Delta_i-r)\cap [0,2\delta-1]\}\cup \{\infty\}$ must be the set of all orders of elements in $\mathcal{M}_{r,i}$,  the coefficients  may satisfy some conditions to keep it. 
We denote by $H$ the set of all such conditions. 
Put $\Gamma(\mathcal{M}_{r,i}):=\{(\Delta_i-r)\cap [0,2\delta-1]\}\cup \{\infty\}$. 
To determine $H$ for a given $\Gamma$-semi-module $\Delta_i$, we introduce the reduction of two elements in $\mathcal{M}_{r,i}$. 
This is an analogue of S-polynomial for given two polynomials (cf.\cite{Cox}).  
Let $\overline{f}$ be an element of $\mathcal{M}_{r,i}$. 
We express the leading coefficient of $\overline{f}$ 
(resp. the leading term of $\overline{f}$) by {\small LC}$(\overline{f})$ (resp. {\small LT}$(\overline{f})$) with respect to the local order 
$\mathrm{ord} (\tau^0) \succ\mathrm{ord} (\tau^1)\succ \mathrm{ord} (\tau^2)\succ \cdots$. 
Define a multiplication of two elements in $\mathcal{O}$ and $\mathcal{M}_{r,i}$ by $\mathcal{O}\times\mathcal{M}_{r,i}\rightarrow \mathcal{M}_{r,i}$, $(f,\overline{h})\rightarrow \overline{fh}$. 
Take $\overline{h}_1,\overline{h}_2\in \mathcal{M}_{r,i}$.  
Let $(\lambda_1,\lambda_2)$ be the element of  $\{(\gamma_1,\gamma_2)\in \Gamma\times \Gamma|\, \gamma_1\cdot\mathrm{ord}(\overline{h}_1)=\gamma_2\cdot\mathrm{ord}(\overline{h}_2)\}$  
that makes the value $\gamma_1\cdot\mathrm{ord}(\overline{h}_1)=\gamma_2\cdot\mathrm{ord}(\overline{h}_2)$ minimal. 
It is clear that such pair is uniquely determined. 
We define the \emph{reduction} of $\overline{h}_1$ and $\overline{h}_2$ by
\begin{equation*}
\mathrm{Red}(\overline{h}_1,\overline{h}_2):=\text{{\small LC}} (\overline{h}_2)\cdot \tau^{\lambda_1}\cdot\overline{h}_1-\text{{\small LC}}(\overline{h}_1) \cdot \tau^{\lambda_2}\cdot\overline{h}_2.
\end{equation*}

\begin{prp}\label{alg3}
The condition set $H$ for  $\mathcal{M}_{r,i}=[ \overline{g}_1,\cdots,\overline{g}_\delta] _k$  is given by the following algorithm$:$

\noindent
\rm{\texttt{INPUT:}} $\{\overline{g}_1,\cdots,\overline{g}_\delta\}$\\
\rm{\texttt{OUTPUT:}} $H$\\
\rm{\texttt{DEFINE:}} $H:=\emptyset$\\
\rm{\texttt{FOR}} each $i,j$ in $\{1,\ldots,2\delta\}$ with $i\neq j$ \rm{\texttt{DO}}\\
\quad $R:=$Red$(\overline{g}_i,\overline{g}_j)$\\
\qquad\rm{\texttt{WHILE}} $\mathrm{ord} (R)<2\delta$ \rm{\texttt{DO}}\\
\qquad\quad\rm{\texttt{IF}} ord$(R)\notin \Gamma(\mathcal{M}_{r,i})$ \rm{\texttt{THEN}} $R:=R-{\mathrm{\small LT}}(R)$ \rm{\texttt{and}}
$H:=H\cup\{\mathrm{{\small LC}}(R)=0\}$\\
\qquad\quad\rm{\texttt{ELSE}} 
$R:=\mathrm{Red}\displaystyle{\Bigg(R,\sum_{\overline{g}_i\in L}t^{\gamma_i}\overline{g}_i\Bigg)}$ for $L=\{\overline{g}_i|\,\exists \gamma_i\in\Gamma\textit{ s.t. }\gamma_i+\textrm{ord}(\overline{g}_i)=\mathrm{ord}(R)\}$
\end{prp}
\begin{proof}
For two distinct basis $\overline{g}_i$ and $\overline{g}_j$, we first compute $R_1:=\mathrm{Red}(\overline{g}_i,\overline{g}_j)$. 
Note that {\small LC}$(R_1)$ is a polynomial with respect to the coefficients in $\overline{g}_i$ and $\overline{g}_j$. 
If $\mathrm{ord}(R_1)\notin \Gamma(\mathcal{M}_{r,i})$, then we must have {\small LC}$(R_1)=0$. 
We add this equation to $H$ and put $R_2:=R_1-${\small LT}$(R_1)$. 
On the other hand, if $\mathrm{ord}(R_1)\in \Gamma(\mathcal{M}_{r,i})$, 
then, for $L_{1}:=\{\overline{g}_i |\,\exists \gamma_i\in\Gamma\textit{ s.t. }\gamma_i+\textrm{ord}(\overline{g}_i)=\mathrm{ord}(R_1)\}$, 
consider the reduction  $R_2:=\mathrm{Red}\left(R_1,\sum_{\overline{g}_i\in L_1}\tau^{\gamma_i}\overline{g}_i\right)$. 
We check whether $\mathrm{ord}(R_2)$ belongs to $\Gamma(\mathcal{M}_{r,i})$ or not. 
Continuing these procedures  successively, we obtain $\mathrm{ord}(R_1)< \mathrm{ord}(R_2)<\cdots$. 
So there exists $s$ which satisfies $\mathrm{ord}(R_{s-1})<2\delta$ and $\mathrm{ord}(R_s)\ge 2\delta$. 
Namely, our procedures terminate in finite steps and we obtain the condition set $H$. 
\end{proof}

We denote by $\mathfrak{m}$ the maximal ideal of $\mathcal{O}$. 
The following theorem follows from Lemma\,\ref{normal form}, \ref{alg1}, Proposition\,\ref{alg2} and \ref{alg3}:

\begin{thm}[Computational algorithm for an affine cell decomposition of $\mathcal{M}_r$]\label{Computational steps}
For a given codimenson $r$, we obtain all affine cells in the decomposition $(\ref{decomposition 2})$ of $\mathcal{M}_r$ by the following finite steps.\\
\rm{\textbf{Step\,1:}} Set $\mathfrak{D}_1=\{\Delta_{1,1}\}$ where $\Delta_{1,1}=\Gamma(\mathfrak{m})$ and find generators of $\Gamma(\mathfrak{m})$ by Lemma\,\ref{alg1}.  \\
\rm{\textbf{Step\,\mathversion{bold}$i$ $(i=2,\ldots,r)$:}} Compute $\mathfrak{D}_{i}$ from $\mathfrak{D}_{i-1}$ by Proposition\,\ref{alg2} 
and, applying Lemma\,\ref{alg1} to each elements in $\mathfrak{D}_{i}$, determine their sets of minimal generators.\\
\rm{\textbf{Step\,r+1:}} For each $\Delta_{r,i}$ in $\mathfrak{D}_r$, determine $\mathcal{I}(\Delta_{r,i})$ by Lemma\,\ref{normal form}. \\
\rm{\textbf{Step\,r+2:}} For each $\mathcal{I}(\Delta_{r,i})$,  compute the condition set $H$ of $\mathcal{M}_{r,i}$ by Proposition\,\ref{alg3} 
and determine $\mathcal{M}_{r,i}$ by Lemma\,\ref{generators for k-vec.}  and $H$. 
\end{thm}

We will use the theory of Gr\"obner basis in the proofs of Theorem\,\ref{main 1} and \ref{main 2}. 
For an ideal $I\subset k[x_1,\ldots,x_n]$,
the ideal $I_m=I\cap k[x_{m+1},\ldots,x_n]$ is called the $m$th \emph{elimination ideal} of $I$. 
The following  theorems are useful:
\begin{thm}[\cite{Cox}, The Elimination Theorem]\label{The Elimination Theorem}
Let $I\in k[x_1,\ldots,x_n]$ be an ideal and let $G$ be a Gr\"obner basis of $I$ with respect to lexicographic ordering where
$x_1>x_2>\cdots>x_n$. 
Then, for every $0\le m\le n$, the set $G_m=G\cap k[x_{m+1},\ldots, x_n]$
is a Gr\"obner basis of the $m$th elimination ideal $I_m$. 
\end{thm}
\begin{thm}[\cite{Cox}, Polynomial Implicitization]\label{Polynomial Implicitization}
Consider the map $F:k^m\rightarrow k^n$ determined by $x_i=f_i(t_1,\ldots,t_m) (i=1,\ldots,n)$. 
Let $I$ be the ideal $(x_1-f_1,\ldots,x_n-f_n)\subset k[t_1,\ldots,t_m,x_1, \ldots, x_n]$ and let $I_m=I\cap k[x_1,\cdots,x_n]$ be the $m$th elimination ideal.
Then the  zero set of $I_m$ is the smallest variety in $k^n$ containing $F(k^m)$. 
\end{thm}

\section{The singularity of type $E_6$}\label{E6}

We prove Theorem\,\ref{main 1} in this section.

\noindent
\textbf{Proof of Theorem\,\ref{main 1}}
Let $\mathcal{O}$ be the local ring $k[[t^3,t^4]]$ of the singularity of type $E_{6}$. 
So we have $\Gamma=\{0,3,4,6,7,8,9,\cdots\}$.
It follows that $\delta=3$ and $c=2\delta=6$. 
For each codimension $r$ $(1\le r\le 6)$, we first determine all components appear in the decomposition\,(\ref{decomposition 2}). 
Performing Step\,1 to Step\,$r$ $(1\le r\le 6)$ in Theorem\ \ref{Computational steps}, we obtain the following datum:
\begin{center}
\begin{tabular}{c|l}
$r$&Elements of $\mathfrak{D}_{r}$\\
\hline
1&$\Delta_{1,1}=\langle 3,4\rangle_\Gamma$\\
\hline
2&$\Delta_{2,1}=\langle 4,6\rangle_\Gamma$, $\Delta_{2,2}=\langle 3,8\rangle_\Gamma$\\
\hline
3&$\Delta_{3,1}=\langle 6,7,8\rangle_\Gamma$, $\Delta_{3,2}=\langle 4,9\rangle_\Gamma$, $\Delta_{3,3}=\langle 3\rangle_\Gamma$\\
\hline
4&$\Delta_{4,1}=\langle 7,8,9\rangle_\Gamma$, $\Delta_{4,2}=\langle 6,8\rangle_\Gamma$, $\Delta_{4,3}=\langle 6,7\rangle_\Gamma$, $\Delta_{4,4}=\langle 4\rangle_\Gamma$\\
\hline
5&$\Delta_{5,1}=\langle 8,9,10\rangle_\Gamma$, $\Delta_{5,2}=\langle 7,9\rangle_\Gamma$, $\Delta_{5,3}=\langle 7,8\rangle_\Gamma$, $\Delta_{5,4}=\langle 6,11\rangle_\Gamma$\\
\hline
6&$\Delta_{6,1}=\langle 9,10,11\rangle_\Gamma$, $\Delta_{6,2}=\langle 8,10\rangle_\Gamma$, $\Delta_{6,3}=\langle 8,9\rangle_\Gamma$, $\Delta_{6,4}=\langle 7,12\rangle_\Gamma$, $\Delta_{6,5}=\langle 6\rangle_\Gamma$\\
\multicolumn{2}{c}{\textbf {Table\,}{\mathversion{bold}$3$.}}
\end{tabular} 
\end{center}

Furthermore, Step\,$r+1$ and Step\,$r+2$ for $r=1,\ldots,6$ yield the following table:

\begin{center}
\begin{tabular}{c|l}
$r$&Components of $\mathcal{M}_r$\\
\hline
1&$\mathcal{M}_{1,1}=[ \tau^2,\tau^3,\tau^5 ]_k$\\
\hline
2&$\mathcal{M}_{2,1}=[ \tau^2,\tau^4,\tau^5 ]_k$, $\mathcal{M}_{2,2}=[ \tau+a\tau^2,\tau^4,\tau^5 ]_k$\\
\hline 
3&$\mathcal{M}_{3,1}=[ \tau^3,\tau^4,\tau^5 ]_k$, $\mathcal{M}_{3,2}=[ \tau+a\tau^3,\tau^4,\tau^5 ]_k$\\
&$\mathcal{M}_{3,3}=[ 1+a\tau+b\tau^5,\tau^3+a\tau^4,\tau^4+a\tau^5 ]_k$  \\
\hline
4&$\mathcal{M}_{4,1}=[ \tau^3,\tau^4,\tau^5 ]_k$, $\mathcal{M}_{4,2}=[ \tau^2+a\tau^3,\tau^4,\tau^5 ]_k$\\
& $\mathcal{M}_{4,3}=[ \tau^2+a\tau^4,\tau^3+b\tau^4,\tau^5 ]_k$, $\mathcal{M}_{4,4}=[ 1+a\tau^2+b\tau^5,\tau^3+a\tau^5,\tau^4 ]_k$\\
\hline
5&$\mathcal{M}_{5,1}=[ \tau^3,\tau^4,\tau^5 ]_k$, $\mathcal{M}_{5,2}=[ \tau^2+a\tau^3,\tau^4,\tau^5 ]_k$\\
&$\mathcal{M}_{5,3}=[ \tau^2+a\tau^4,\tau^3+b\tau^4,\tau^5 ]_k$, $\mathcal{M}_{5,4}=[ \tau+a\tau^2+b\tau^3,\tau^4,\tau^5 ]_k$\\
\hline
6&$\mathcal{M}_{6,1}=[ \tau^3,\tau^4,\tau^5 ]_k$, $\mathcal{M}_{6,2}=[ \tau^2+a\tau^3,\tau^4,\tau^5 ]_k$, $\mathcal{M}_{6,3}=[ \tau^2+a\tau^4,\tau^3+b\tau^4,\tau^5 ]_k$\\
& $\mathcal{M}_{6,4}=[ \tau+a\tau^2+b\tau^3,\tau^4,\tau^5 ]_k$, $\mathcal{M}_{6,5}=[ 1+a\tau+b\tau^2+c\tau^5,\tau^3+(b-a^2)\tau^5,\tau^4+a\tau^5 ]_k$\\
\multicolumn{2}{c}{\textbf {Table\,}{\mathversion{bold}$4$.}}
\end{tabular}
\end{center}
In Table\,4,   $a,b,c\in k$. 
Here we only explain the case of  $\mathcal{M}_6$. 
This is the most complicated case for the singularity of type $E_6$. 
The other cases can be teated in the similar manner. 
Since all element in $\mathcal{M}_{6,i}$ has generators of the form (\ref{basis}), their represent matrices have same form. 
So we just express them by $A_i$. 
The matrices $A_i$ $(i=1\,\ldots,5)$ are calculated from Table\,4  as follows:
\begin{align*}
A_1&=
\begin{pmatrix}
0&0&0&1&0&0\\
0&0&0&0&1&0\\
0&0&0&0&0&1\\
\end{pmatrix},\ 
&A_2=
\begin{pmatrix}
0&0&1&a&0&0\\
0&0&0&0&1&0\\
0&0&0&0&0&1\\
\end{pmatrix},\ 
&A_3=
\begin{pmatrix}
0&0&1&0&a&0\\
0&0&0&1&b&0\\
0&0&0&0&0&1\\
\end{pmatrix}\\ 
A_4&=
\begin{pmatrix}
0&1&a&b&0&0\\
0&0&0&0&1&0\\
0&0&0&0&0&1\\
\end{pmatrix},\ 
&A_5=
\begin{pmatrix}
1&a&b&0&0&c\\
0&0&0&1&0&b-a^2\\
0&0&0&0&1&a
\end{pmatrix}\ &
\end{align*}
The Pl\"ucker coordinates $\pi_{ijk}$ for each $\mathcal{M}_{6,i}$ are defined to be the determinants which consists of $i,j$ and $k$th columns of $A_i$ $(1\le i<j<k\le 6)$.  
They are calculated   as follows:
\begin{align*}
\mathcal{M}_{6,1}:&\pi_{456}=1,\,\pi_{ijk}=0\ \text{for $(i,j,k)\neq(4,5,6)$}\\
\mathcal{M}_{6,2}:&\pi_{356}=1,\,\pi_{456}=a,\,\pi_{ijk}=0\ \text{for $(i,j,k)\neq(3,5,6),\,(4,5,6)$}\\
\mathcal{M}_{6,3}:&\pi_{346}=1,\,\pi_{356}=b,\,\pi_{456}=-a, \pi_{ijk}=0\ \text{for $(i,j,k)\neq(3,4,6),\,(3,5,6),\,(4,5,6)$}\\
\mathcal{M}_{6,4}:&\pi_{256}=1,\,\pi_{356}=a,\,\pi_{456}=b, \pi_{ijk}=0\ \text{for $(i,j,k)\neq(2,5,6),\,(3,5,6),\,(4,5,6)$}\\
\mathcal{M}_{6,5}:&\pi_{145}=1,\,\pi_{146}=a,\,\pi_{156}=a^2-b,\,\pi_{245}=a,\,\pi_{246}=a^2,\,\pi_{256}=a^3-ab,\,\pi_{345}=a\\
&\pi_{346}=a^2,\,\pi_{345}=a,\,\pi_{346}=a^2,\,\pi_{356}=a^2b-b^2,\,\pi_{456}=c,\,\pi_{ijk}=0\ \text{for the others}
\end{align*}
By using these Pl\"ucker coordinates, we can check that
\begin{align}\label{union}
&\mathcal{M}_{6,1}\cup\mathcal{M}_{6,2}\cup\mathcal{M}_{6,3}\cong\mathcal{M}_{6,1}\cup\mathcal{M}_{6,2}\cup\mathcal{M}_{6,4}\cong\mathbb{P}^2,\\\label{intersection}
&(\mathcal{M}_{6,1}\cup\mathcal{M}_{6,2}\cup\mathcal{M}_{6,3})\cap(\mathcal{M}_{6,1}\cup\mathcal{M}_{6,2}\cup\mathcal{M}_{6,4})=\mathcal{M}_{6,1}\cup\mathcal{M}_{6,2}\cong\mathbb{P}^1.
\end{align}
We  calculate  the defining equations of  $\mathcal{M}_{6,5}$ to show  $\mathcal{M}_6=\overline{\mathcal{M}_{6,5}}$. 
Let $I$ be an ideal which is generated by the following polynomials in $k[\pi_{ijk}|\,1\le i<j<k\le 6,\,(i,j,k)\neq(1,4,5)]$:
\begin{align*}
&\pi_{146}-a,\,\pi_{156}-a^2+b,\,\pi_{245}-a,\,\pi_{246}-a^2,\,\pi_{256}-a^3+ab,\,\pi_{345}-a,\\
&\pi_{346}-a^2,\,\pi_{345}-a,\,\pi_{346}-a^2,\,\pi_{356}-a^2b+b^2,\,\pi_{456}-c
\end{align*}
By Pl\"ucker coordinates of $\mathcal{M}_{6,5}$ listed above and Theorem\,\ref{Polynomial Implicitization}, $\mathcal{M}_{6,5}$ is defined by the third elimination ideal $I_3$ of $I$. 
We compute the Gr\"obner basis of $I$ with respect to a lexicographic ordering where $a\succ b\succ c\succ \pi_{123}\succ \cdots\succ \pi_{456}$. 
By Theorem\,\ref{The Elimination Theorem}, the elements of the Gr\"obner basis not involving $a,b,c$ form a basis of $I_3$. 
This computation was done by the computer algebra system ^^ ^^ Singular'' 
(see \cite{GP} for the usage of Singular). 
Furthermore, homogenizing the basis of $I_3$ at $\pi_{145}$, we obtain the defining equations of a projective three fold named $X_2$ as follows:
\begin{align*}
&\pi_{345}^3+\pi_{145}\pi_{345}\pi_{356}-\pi_{145}\pi_{346}^2=0,\,\pi_{256}\pi_{345}^2-\pi_{145}\pi_{346}\pi_{356}=0,\\
&\pi_{145}\pi_{256}\pi_{346}^3-\pi_{345}^3\pi_{356}-\pi_{345}^3\pi_{346}^2-2\pi_{145}\pi_{345}^2\pi_{356}^2\\
&\qquad\qquad\qquad\qquad\qquad\qquad\qquad\qquad-\pi_{145}\pi_{345}\pi_{346}^2\pi_{356}+\pi_{145}\pi_{346}^4-\pi_{145}^2\pi_{356}^3=0,\\
&\pi_{145}\pi_{256}\pi_{345}\pi_{356}-\pi_{145}\pi_{256}\pi_{346}^2+\pi_{345}^3\pi_{346}+2\pi_{145}\pi_{345}\pi_{346}\pi_{356}-\pi_{145}\pi_{346}^3=0,\\
&\pi_{145}\pi_{256}\pi_{345}\pi_{346}-\pi_{345}^4-2\pi_{145}\pi_{345}^2\pi_{356}+\pi_{145}\pi_{345}\pi_{346}^2-\pi_{145}^2\pi_{356}^2=0,\\
&\pi_{145}\pi_{246}\pi_{356}-\pi_{145}\pi_{256}\pi_{346}+\pi_{345}^3+\pi_{145}\pi_{345}\pi_{356}-\pi_{145}\pi_{346}^2=0,\\
&\pi_{246}\pi_{346}-\pi_{256}\pi_{345}-\pi_{345}\pi_{346}=0,\,\pi_{246}\pi_{345}-\pi_{345}^2-\pi_{145}\pi_{356}=0,\\
&\pi_{246}^3-\pi_{246}^2\pi_{345}-\pi_{145}\pi_{256}^2-\pi_{145}\pi_{256}\pi_{346}=0,\pi_{245}\pi_{356}-\pi_{256}\pi_{345}=0,\\
&\pi_{245}\pi_{346}-\pi_{246}\pi_{345}=0,\,\pi_{245}\pi_{345}-\pi_{145}\pi_{346}=0,\,\pi_{245}\pi_{256}+\pi_{245}\pi_{346}-\pi_{246}^2=0,\\
&\pi_{245}\pi_{246}-\pi_{245}\pi_{345}-\pi_{145}\pi_{256}=0,\,\pi_{245}^2-\pi_{145}\pi_{246}=0,\\
&\pi_{145}\pi_{156}-\pi_{245}^2+\pi_{145}\pi_{345}=0,\,\pi_{146}-\pi_{245}=0
\end{align*}

By using computer algebra system ^^ ^^ Maple'', we check that  $\mathcal{M}_6\setminus\mathcal{M}_{6,5}$ is defined by these equations with $\pi_{145}=0$. 
This fact implies $\mathcal{M}_6=\overline{\mathcal{M}}_{6,5}=X_2$. 
So the variety $X_2$ is irreducible. 
The dimensions $\dim \mathcal{M}_6=\dim \mathcal{M}_{6,5}=3$ and the rationality of $\mathcal{M}_6$ also follow from Proposition\,\ref{affine space} and Corollary\,\ref{rationality} respectively. 
We also conclude that Sing($\mathcal{M}_6)=\mathbb{P}^1$ by ($\ref{union}$) and ($\ref{intersection}$). 

We add some comments for the other cases. 
One can check that $\mathcal{M}_4$ consists of two components, one is $\mathbb{P}^1$ and the other is a rational surface 
$X_1$ given by $\pi_{156}+\pi_{345}=0$ and $\pi_{145}\pi_{356}+\pi_{345}^2=0$. 
We also see that $\mathcal{M}_4$, $\mathcal{M}_5$ and $\mathcal{M}_6$ possess same $\mathbb{P}^1$ in common, as their singular locus. 
$\square$

The irreducibility of $\mathcal{M}_6$  can  be proven by the known results of the compactified Jacobians. 
In next section, we use them to show the irreducibility of PS variety for the singularity of type $E_8$.

\section{The singularity of type $E_8$}\label{E8}

Consider the curve singularity of type $E_8$ in this section.
In order to prove Theorem\,\ref{main 2}, we recall some results about compactified Jacobian $\overline{JC}$ for a singular complete algebraic curve $C$. 
The compactified Jacobian $\overline{JC}$ was defined by Rego. 
\begin{Def}[\cite{R}]
The compactified Jacobian $\overline{JC}$ of $C$ consists of all torsion free sheaves $\mathcal{F}$ of rank $1$ and degree $0$ on $C$
$($i.e. $\chi (\mathcal{F})=1-g_a(C))$. 
\end{Def}

\noindent
The following facts about compactified Jacobians are known: 
\begin{thm}[\cite{AIK}, \cite{R}]\label{irreducibility of the compactified Jacobian }
The compactified Jacobian $\overline{JC}$ is irreducible if and only if $\mathrm{Sing} (C)$ consits of plane curve singularities. 
\end{thm}
\begin{thm}[\cite{B}]\label{B}
For a rational unibranched curve $C$, its compactified Jacobian  is homeomorphic to the direct product of compact spaces, 
the Jacobi factors $\overline{JC}_p$ where  $p\in\mathrm{Sing} (C)$. 
\end{thm}

\noindent
The Jacobi factor for a curve singularity was introduced by Rego in \cite{R} (see also \cite{P}). 
It was defined to be the punctual Hilbert scheme of degree $2\delta$ for the singularity. 
\begin{rem}\label{rem}
For the case $r=2\delta$, Corollary\,\ref{decomposition cor} gives an affine cell decomposition of the Jacobi factor of a given singularity. 
So the Euler number of Jacobi factor is given by $l_{2\delta}$ in $(\ref{decomposition 2})$.  
\end{rem}

\noindent
\textbf{Proof of Theorem\,\ref{main 2}} 
Let $\mathcal{O}$ be the ring $k[[t^3,t^5]]$. 
We have $\Gamma=\{0,3,5,6,8,9,\cdots\}$.
It follows that $\delta=4$ and $c=2\delta=8$. 
By Theorem\ \ref{Computational steps}, we obtain the following two tables:
\begin{center}
\begin{tabular}{c|l}
$r$&Elements of $\mathfrak{D}_{i}$\\
\hline
1&$\Delta_{1,1}=\langle 3,5\rangle_\Gamma$\\
\hline
2&$\Delta_{2,1}=\langle 5,6\rangle_\Gamma$, $\Delta_{2,2}=\langle 3,10\rangle_\Gamma$\\
\hline
3&$\Delta_{3,1}=\langle 6,8,10\rangle_\Gamma$, $\Delta_{3,2}=\langle 5,9\rangle_\Gamma$, $\Delta_{3,3}=\langle 3\rangle_\Gamma$\\
\hline
4&$\Delta_{4,1}=\langle 8,9,10\rangle_\Gamma$, $\Delta_{4,2}=\langle 6,10\rangle_\Gamma$, $\Delta_{4,3}=\langle 6,8\rangle_\Gamma$, $\Delta_{4,4}=\langle 5,12\rangle_\Gamma$\\
\hline
5&$\Delta_{5,1}=\langle 9,10,11\rangle_\Gamma$, $\Delta_{5,2}=\langle 8,10,12\rangle_\Gamma$, $\Delta_{5,3}=\langle 8,9\rangle_\Gamma$, $\Delta_{5,4}=\langle 6,13\rangle_\Gamma$, $\Delta_{5,5}=\langle 5\rangle_\Gamma$\\
\hline
6&$\Delta_{6,1}=\langle 10,11,12\rangle_\Gamma$, $\Delta_{6,2}=\langle 9,11,13\rangle_\Gamma$, $\Delta_{6,3}=\langle 9,10\rangle_\Gamma$, $\Delta_{6,4}=\langle 8,12\rangle_\Gamma$\\
&$\Delta_{6,5}=\langle 8,10\rangle_\Gamma$, $\Delta_{6,6}=\langle 6\rangle_\Gamma$\\
\hline
7&$\Delta_{7,1}=\langle 11,12,13\rangle_\Gamma$, $\Delta_{7,2}=\langle 10,12,14\rangle_\Gamma$, $\Delta_{7,3}=\langle 10,11\rangle_\Gamma$, $\Delta_{7,4}=\langle 9,13\rangle_\Gamma$\\
&$\Delta_{7,5}=\langle 9,11\rangle_\Gamma$, $\Delta_{7,6}=\langle 8,15\rangle_\Gamma$\\
\hline
8&$\Delta_{8,1}=\langle 12,13,14\rangle_\Gamma$, $\Delta_{8,2}=\langle 11,13,15\rangle_\Gamma$, $\Delta_{8,3}=\langle 11,12\rangle_\Gamma$, $\Delta_{8,4}=\langle 10,14\rangle_\Gamma$\\
&$\Delta_{8,5}=\langle 10,12\rangle_\Gamma$, $\Delta_{8,6}=\langle 9,16\rangle_\Gamma$,  $\Delta_{8,7}=\langle 8\rangle_\Gamma$\\
\multicolumn{2}{c}{\textbf {Table\,}{\mathversion{bold}$5$.}}
\end{tabular} 
\end{center}

\begin{center}
\begin{tabular}{c|l}
$r$&Components of $\mathcal{M}_r$\\
\hline
1&$\mathcal{M}_{1,1}=[ \tau^2,\tau^4,\tau^5,\tau^7 ]_k$\\
\hline
2&$\mathcal{M}_{2,1}=[ \tau^3,\tau^4,\tau^6,\tau^7 ]_k$, $\mathcal{M}_{2,2}=[ \tau+a\tau^3,\tau^4,\tau^6,\tau^7 ]_k$\\
\hline 
3&$\mathcal{M}_{3,1}=[ \tau^3,\tau^5,\tau^6,\tau^7 ]_k$, $\mathcal{M}_{3,2}=[ \tau^2+a\tau^3,\tau^5,\tau^6,\tau^7 ]_k$\\
&$\mathcal{M}_{3,3}=[ 1+a\tau^2+b\tau^7,\tau^3-a^2\tau^7,\tau^5+a\tau^7,\tau^6 ]_k$  \\
\hline
4&$\mathcal{M}_{4,1}=[ \tau^4,\tau^5,\tau^6,\tau^7 ]_k$, $\mathcal{M}_{4,2}=[ \tau^2+a\tau^4,\tau^5,\tau^6,\tau^7 ]_k$\\
& $\mathcal{M}_{4,3}=[ \tau^2+a\tau^6,\tau^4+b\tau^6,\tau^5,\tau^7 ]_k$, $\mathcal{M}_{4,4}=[ \tau+a\tau^2+b\tau^5,\tau^4+a\tau^5,\tau^6,\tau^7 ]_k$\\
\hline
5&$\mathcal{M}_{5,1}=[ \tau^4,\tau^5,\tau^6,\tau^7 ]_k$, $\mathcal{M}_{5,2}=[ \tau^3+a\tau^4,\tau^5,\tau^6,\tau^7 ]_k$\\
&$\mathcal{M}_{5,3}=[ \tau^3+a\tau^5,\tau^4+b\tau^5,\tau^6,\tau^7 ]_k$, $\mathcal{M}_{5,4}=[ \tau+a\tau^3+b\tau^5,\tau^4,\tau^6,\tau^7 ]_k$\\
& $\mathcal{M}_{5,5}=[ 1+a\tau+b\tau^4+c\tau^7,\tau^3+a\tau^4+b\tau^7,\tau^5-a^2\tau^7,\tau^6+a\tau^7 ]_k$\\
\hline
6&$\mathcal{M}_{6,1}=[ \tau^4,\tau^5,\tau^6\tau^7 ]_k$, $\mathcal{M}_{6,2}=[ \tau^3+a\tau^4,\tau^5,\tau^6,\tau^7 ]_k$\\ &$\mathcal{M}_{6,3}=[ \tau^3+a\tau^5,\tau^4+b\tau^5,\tau^6,\tau^7 ]_k$, $\mathcal{M}_{6,4}=[ \tau^2+a\tau^3+b\tau^4,\tau^5,\tau^6,\tau^7 ]_k$\\
&$\mathcal{M}_{6,5}=[ \tau^2+a\tau^3+b\tau^6,\tau^4+c\tau^6,\tau^5+a\tau^6,\tau^7 ]_k$\\
&$\mathcal{M}_{6,6}=[ 1+a\tau^2+b\tau^4+c\tau^7,\tau^3+(b-a^2)\tau^7,\tau^5+a\tau^7,\tau^6 ]_k$\\
\hline
7&$\mathcal{M}_{7,1}=[ \tau^4,\tau^5,\tau^6\tau^7 ]_k$, $\mathcal{M}_{7,2}=[ \tau^3+a\tau^4,\tau^5,\tau^6,\tau^7 ]_k$\\ &$\mathcal{M}_{7,3}=[ \tau^3+a\tau^5,\tau^4+b\tau^5,\tau^6,\tau^7 ]_k$, 
$\mathcal{M}_{7,4}=[ \tau^2+a\tau^3+b\tau^4,\tau^5,\tau^6,\tau^7 ]_k$\\
&$\mathcal{M}_{7,5}=[ \tau^2+a\tau^3+b\tau^6,\tau^4+c\tau^6,\tau^5+a\tau^6,\tau^7 ]_k$\\
&$\mathcal{M}_{7,6}=[ \tau+a\tau^2+b\tau^3+c\tau^5,\tau^4+a\tau^5,\tau^6,\tau^7 ]_k$\\
\hline
8&$\mathcal{M}_{8,1}=[ \tau^4,\tau^5,\tau^6\tau^7 ]_k$, $\mathcal{M}_{8,2}=[ \tau^3+a\tau^4,\tau^5,\tau^6,\tau^7 ]_k$\\ 
&$\mathcal{M}_{8,3}=[ \tau^3+a\tau^5,\tau^4+b\tau^5,\tau^6,\tau^7 ]_k$, $\mathcal{M}_{8,4}=[ \tau^2+a\tau^3+b\tau^4,\tau^5,\tau^6,\tau^7 ]_k$\\
&$\mathcal{M}_{8,5}=[ \tau^2+a\tau^3+b\tau^6,\tau^4+c\tau^6,\tau^5+a\tau^6,\tau^7 ]_k$\\
&$\mathcal{M}_{8,6}=[ \tau+a\tau^2+b\tau^3+c\tau^5,\tau^4+a\tau^5,\tau^6,\tau^7 ]_k$\\
&$\mathcal{M}_{8,7}=[ 1+a\tau+b\tau^2+c\tau^4+d\tau^7,\tau^3+a\tau^4+(c+a^2b-b^2)\tau^7,\tau^5+(b-a^2)\tau^7, \tau^6+a\tau^7 ]_k$\\
\multicolumn{2}{c}{\textbf {Table\,}{\mathversion{bold}$6$.}}
\end{tabular}
\end{center}
In Table\,6,  we have $a,b,c\in k$. 
The analysis for $\mathcal{M}_i$ for $i=1,\ldots,8$ proceeds similarly as in the proof of Theorem\,\ref{main 1}, except the irreducibility of $\mathcal{M}_8$. 
The defining equations of  $\mathcal{M}_{8,7}$ are  too many to analyze by the using Gr\"obner basis. 
So we use Theorem\,\ref{irreducibility of the compactified Jacobian } and \ref{B} to show  the irreducibility of $\mathcal{M}_8$. 
Let $C$ be a rational curve with the curve singularity of type $E_8$ as its unique singularity. 
Since the compactified Jacobian $\overline{JC}$ is irreducible by Theorem\,\ref{irreducibility of the compactified Jacobian }, 
the irreducibility of  the PS variety $\mathcal{M}_8$ follows from Theorem\,\ref{B}. 
Finally, we conclude that $\mathcal{M}_{8}$ is rational by  Corollary\,\ref{rationality}. 

For the other cases, we only mentioned that the punctual Hilbert schemes $\mathcal{M}_6$, $\mathcal{M}_7$ and $\mathcal{M}_8$ possess same $\mathbb{P}^2\cup\mathbb{P}^2$ in common, as their singular locus.  
$\square$
\begin{rem}
In $\cite{B}$, Beauville proved that the Euler number of the Jacobi factor for the curve singularities of types $E_6$ $($resp. $E_8)$ is $5$ $($resp. $7)$.  
As in  Remark\,$\ref{rem}$, the Euler number of the Jacobi factor of a curve singularity equals the number of its affine cells  $l_{2\delta}=\sharp\mathfrak{D}_{2\delta}$. 
So, for the singularities of types $E_6$ and $E_8$, the Euler numbers of the Jacobi factors are also  derived from Table\,$4$ and $6$.
\end{rem}

\noindent
\small 
Yoshiki S\={o}ma\\
Ch\={o}fu Minami Metroporitan High School\\
6-2-1 Tamagawa\\
Ch\={o}fushi Tokyo 182-0025 Japan.\\

\noindent
\small 
Masahiro Watari\\
Division of General Subjects\\
Tsuyama National College of Technology\\
624 Numa\\
Tsuyamashi Okayama 708-8509  Japan.\\
E-mail:watari@tsuyama-ct.ac.jp

\begin{thebibliography}{9}
\bibitem{AIK}
A. Altman, A. Iarrobino, S. Kleiman, 
Irreducibility of the compactified Jacobian.
Real and Complex Singularities, Proc. Nordic Summer Sch., Symp. Math., Oslo 1976 (1977), 1-12. 
\bibitem{B}
A. Beauville. 
Counting rational curves on $K3$-surfaces. 
Dule Math. J. \textbf{97}(1999), 99-108. 
\bibitem{Cox}
D. Cox,  J. Little,  D. O'Shea,  
IDEALS, VARIETIES, AND ALGORITHMS. An introduction to computational algebraic geometry and commutative algebra. Third edition. 
Undergraduate Texts in Mathematics. Springer, New York (2007).
\bibitem{GP}
G-M. Greuel, G. Pfister, 
A Singular Introduction to Commutative Algebra
Springer-Verlag Berlin Heidelberg New York (2002).
\bibitem{P}
J. Piontkowski, 
Topology of the compactified Jacobians of singular curves.
Math. Z. \textbf{255} (2007), 195-226.  
\bibitem{PS}
G. Pfister, J.H.M. Steenbrink, 
Reduced Hilbert schemes for irreducible curve singularities. 
J. Pure and Applied Algebra. \textbf{77} (1992), 103-116. 
\bibitem{R}
C. Rego, 
The compactified Jacobian. 
Ann. Sci. \'{E}c. Norm. Sup\'{e}r., IV. S\'{e}r. \textbf{13} (1980), 211-223. 
\bibitem{Serre}
J. P. Serre, 
Groupes Alg\'ebriques et Corps de Classes. 
Hermann, Paris (1959).
\bibitem{S}
Y. S\={o}ma,
Hilbert scheme of $r$-points for the simple plane curve singularities, (in Japanese)        
Master Thesis, Saitama university (2010).  
\end{thebibliography}
\end{document}